\documentclass[amscd,amssymb,verbatim,a4paper]{amsart}

\theoremstyle{plain}
\newtheorem{theorem}{Theorem}[section]
\newtheorem{lemma}[theorem]{Lemma}
\newtheorem{proposition}[theorem]{Proposition}
\newtheorem{corollary}[theorem]{Corollary}

\theoremstyle{definition}

\theoremstyle{remark}

\newtheorem{notations}[theorem]{Notations}

\newcommand{\ic}{\ensuremath{\mathcal{I}}}
\newcommand{\oc}{\ensuremath{\mathcal{O}}}

\newcommand{\fc}{\ensuremath{\mathcal{F}}}

\newcommand{\lc}{\ensuremath{\mathcal{L}}}
\newcommand{\Sc}{\ensuremath{\mathcal{S}}}
\newcommand{\tc}{\ensuremath{\mathcal{T}}}
\newcommand{\nc}{\ensuremath{\mathcal{N}}}

\newcommand{\Pt}{\mathbb{P}^3}

\newcommand{\Pq}{\mathbb{P}^4}
\newcommand{\Pcq}{\mathbb{P}^5}

\newcommand{\Pu}{\mathbb{P}^1}
\newcommand{\Pn}{\mathbb{P}^n}

\begin{document}

\title[On subcanonical surfaces of $\Pq$.]{On subcanonical surfaces of $\Pq$.}

\author{Ellia Ph.}
\address{Dipartimento di Matematica, via Machiavelli 35, 44100 Ferrara (Italy)}
\email{phe@dns.unife.it}

\author{Franco D.}
\address{Dipartimento di Matematica e Applicazioni "R. Caccioppoli", Univ. Napoli "Federico II", Ple Tecchio 80, 80125 Napoli (Italy)}
\email{davide.franco@unina.it}

\author{Gruson L.}
\address{D\'epartement de Math\'ematiques, Universit\'e de Versailles-St Quentin, Versailles (France)}

\date{8/07/2004}

\maketitle

\section*{Introduction.}

We work over an algebraically closed field of characteristic zero. It is guessed that there is no indecomposable, unstable rank two vector bundle on $\Pn$, $n \geq 4$. Two of us (\cite{EF}) made an attempt in this direction and we present here a small improvement of it, namely:

\begin{theorem}
\label{thm_subcan_quartic} 
Let $S \subset \Pq$ be a smooth surface with $\omega _S \simeq \oc _S(e)$. If $h^0(\ic _S(4)) \neq 0$,
then $S$ is a complete intersection.  
\end{theorem}

Recall that there exist abelian surfaces of degree ten in $\Pq$ lying on hyperquintics.\\
Here is an outline of the proof. Let $E$ be a rank two vector bundle of which $S$ is a section. The assumption is that $E(-e-1)$ has a section, $Z$, which is of degree $\leq 9$ with $\omega _Z \simeq \oc _Z(-e-2)$; as the cases $e \leq 1$ follow from \cite{BC}, we can also assume that all its components are non reduced. We first show that $Z_{red}$ contains no quadric. We dispose of the case when $Z_{red}$ is a plane $\Pi$ by looking at $E|_{\Pi}$ and deducing that the multiple structure is primitive in codimension one which implies $e \leq 1$. By Lemma \ref{inSing} an irreducible component of degree $\geq 3$ of $Z_{red}$ has to be contained in the singular locus of $\Sigma$ (the hyperquartic containing $S$). It follows that any such a component has degree three. So we are left with the cases where $Z_{red}$ has an irreducible component of degree three. A case by case argument (using results on double structures, see Section \ref{dble-str-surfaces}) concludes the proof.\\
Let us point out some differences between this approach and that of \cite{EF}:
we cannot use, as in \cite{EF}, the linear normality of the hyperplane section;
we have less information on $Z_{red}\cap Sing(\Sigma )$ and, also, no information about $Pic(S)$.


\section{Preliminaries.}

\begin{notations}
\label{gen_setting}
Let $S\subset \Pq$ be a smooth surface of degree $d$ with $\omega _S \simeq \oc _S(e)$ and $h^0(\ic _S(4))\neq 0$. We suppose that $S$ is not a complete intersection. By a result of Koelblen (\cite{K}), we may assume $h^0(\ic _S(3))=0$ (see \cite{E}), hence we may assume that $S \subset \Sigma$ where $\Sigma$ is an integral quartic hypersurface. Also let's us observe that, by \cite{BC}, {\it we may assume $e\geq 2$}, since the only subcanonical surfaces in $\Pq$ with $e \leq 1$ are the abelian surfaces of degree ten, but such surfaces do not lie on a quartic. Since $S$ is subcanonical we may associate to it a rank two vector bundle:
$$0 \to \oc \stackrel{s}{\rightarrow} E \to \ic _S(e+5) \to 0 \:\:(*)$$
here $E$ is a rank two vector bundle with Chern classes: $c_1(E)=e+5$, $c_2(E)=d$. Since $h^0(\ic _S(4))\neq 0$, $h^0(E(-e-1))\neq 0$ and this is the least twist of $E$ having a section. Since $S$ is not a complete intersection, $E(-e-1)$ has a section vanishing in codimension two:
$$0 \to \oc \stackrel{\sigma}{\rightarrow} E(-e-1) \to \ic _Z(-e+3) \to 0 \:\: (**)$$
here $Z$ is a l.c.i. surface of degree $d(Z)=c_2(E(-e-1))=d-4e-4$ and with $\omega _Z \simeq \oc _Z(-e-2)$.\\
We will denote by $Y$ a general hyperplane section of $Z$; $Y$ is a l.c.i. curve with $\omega _Y \simeq \oc _Y(-e-1)$. In particular: $p_a(Y) = \frac{-d(Z)(e+1)}{2}+1$. Let $Y_{red}=Y_1\cup \cdots \cup Y_r$ be the decomposition into irreducible components. Making a primary decomposition we can write: $Y = \tilde{Y}_1\cup \cdots \cup \tilde{Y}_r$, where $\tilde{Y}_i$ is a locally Cohen-Macaulay (and generically l.c.i.) subscheme with support $Y_i$. So $\tilde{Y}_i$ is a locally Cohen-Macaualy multiple structure of multiplicity $m_i$ on the integral curve $Y_i$. The multiplicity $m_i$ is determined by: $deg(\tilde{Y}_i)=m_i.deg(Y_i)$.
\end{notations}

\begin{lemma}
\label{gen_Z}
With notations as in \ref{gen_setting}, we have:\\
(i) If $e \geq 2$, $Z$ is non-reduced, more precisely $m_i\geq 2$, $\forall i$.\\
(ii) $d(Z) \leq 9$; moreover $d(Z) < 9$ if $dim(S \cap Z)>0$.
\end{lemma}

\begin{proof}
(i) For this we argue as in \cite{EF}, Lemma 2.6.\\
(ii) This is \cite{EP} Lemme 1 (see also \cite{EF}, Lemma 2.2). For the last statement following the proof of \cite{EP} Lemme 1, we have $d(Z).d \leq (\sigma -1)^2d-(3H+K)D < (\sigma -1)^2d$ if $D$ is a divisor in $S \cap Jac(\Sigma )$ ($\sigma = d(\Sigma )$).
\end{proof}

\begin{lemma}
\label{inSing}
Let $S\subset \Pq$ be a smooth, non complete intersection surface with $\omega _S\simeq \oc _S(e)$, with $e \geq 2$. Assume $S \subset \Sigma$, where $\Sigma$ is an integral quartic hypersurface. Let $T$ be an irreducible component of $Z_{red}$, if $deg(T) \geq 3$, then $T \subset Sing(\Sigma )$.\\
It follows in particular that every irreducible component of $Z_{red}$ has degree at most three.
\end{lemma}

\begin{proof}
We argue exactly as in \cite{EF} Prop. 3.6.\\
The last statement follows from the fact that an irreducible quartic plane curve has at most three singular points.
\end{proof}

\label{deg2}

\begin{lemma}
\label{E_H}
For any hyperplane $H\subset \Pq $ we have $h^0(E_H(-e-3))=0$. In particular if $dim(Z\cap H)=2$ then the bidimensional component of $Z\cap H$ is a simple plane.
\end{lemma}
\begin{proof}
Since $h^1(E(-e-4)) = h^1(\ic_X(1))=0=h^0(E(-e-3))$, the second equality coming from Severi's Theorem, the sequence
$$
H^0(E(-e-3)) \to H^0(E_H(-e-3)) \to H^1(E(-e-4))
$$
gives the statement.
\end{proof}

\begin{corollary}
\label{irr-comp-degTwo}
With the usual notations, $Z_{red}$ has no irreducible component of degree two, more precisely $Z_{red}$ doesn't contain any  (not necessarily irreducible) degenerate surface of degree $\geq 2$.
\end{corollary}

\begin{lemma}
\label{Z=Pi}
If every irreducible component of $Z_{red}$ has degree one, then $Z_{red}$ is a plane.
\end{lemma}

\begin{proof}
Assume $Z_{red}=\Pi _1 \cup \cdots \cup \Pi _r$, $r>1$, where the $\Pi _i$'s are planes. By Lemma \ref{irr-comp-degTwo}, for every $(i,j)$ with $i \neq j$, we have $dim(\Pi _i \cap \Pi _j)=0$. By \cite{Ha3} it follows that $Z$ can't be locally Cohen-Macaulay. This is absurd since $Z$ is locally a complete intersection.
\end{proof}

By \ref{inSing}, \ref{irr-comp-degTwo} and \ref{Z=Pi}  we are left the following possibilities:
\begin{enumerate}
\item $Z_{red}$ is a plane $\Pi$.
\item $Z_{red}$ is a skew cubic $T$.
\item $Z_{red}$ is the union of a plane with a skew cubic.
\end{enumerate}


\section{ $Z_{red}$ is not a plane.}

In this section we assume that $\Pi \simeq Z_{red}$ is a plane. We start with a further consequence of \ref{E_H}:

\begin{corollary}
\label{noinfn}
If $Z_{red}$ contains a plane $\Pi $ then  $Z$ does not contain $\Pi ^{(1)}$, the infinitesimal neighbourhood of $\Pi$.
\end{corollary}
\begin{proof}
If $Z$ would contain the infinitesimal neighbourhood of $Z$ then any hyperplane $H\supset \Pi $ would contain the double of $\Pi $ contradicting \ref{E_H}.
\end{proof}

Denote by $Z_2\subset Z$ the, unique by \ref{noinfn}, double structure supported on $\Pi $. 

The scheme $Z_2$ comes, via Ferrand's construction, from a surjection
$\nc _{\Pi}^* \to \Sc $ with $\Sc $  a torsionfree sheaf of rank one (on $\Pi$). In the sequel we will show that $Z$ is primitive (i.e. locally contained in a smooth hypersurface) outside of a zero-dimensional scheme; this will force $e=1$.

\begin{lemma}
\label{c1sc}
The first Chern class of $\Sc$ is $\geq 0$.
\end{lemma}
\begin{proof}
Since $\Sc $ is a coker of $\nc _{\Pi}^* \simeq \oc _{\Pi}(-1)\oplus \oc _{\Pi}(-1)$ then $c_1(\Sc )\geq -1$.
If $c_1(\Sc )= -1$ then $Z_2$ would be the double of $\Pi $ inside some hyperplane $H\supset \Pi $ contradicting
\ref{E_H}
\end{proof}

Consider the natural map $\alpha : E_{\Pi}(-4) \simeq \nc^*_Z\mid_{\Pi}\to \nc ^* _{\Pi}$ and set $\tc :=Im(\alpha)$. Since $\Pi ^{(1)}\not \subset Z$ then $\tc $ is a torsionfree sheaf of rank one on $\Pi$.
Moreover, by \ref{c1sc}, we have
$$
0 \to \oc_{\Pi}(-2-l) \to \nc ^*_{\Pi}\to \Sc \to 0
$$
with $l\geq 0$ and $\beta: \tc\hookrightarrow \oc(-2-l)$.

\begin{lemma}
\label{l=0} We have $l=0$ and  $\beta $ is an isomorphism outside from a zero dimensional subscheme, $\Phi$, of $\Pi $. 
\end{lemma}
\begin{proof} 
We have
$$
0 \to \oc_{\Pi}(e-1-k) \to E_{\Pi}(-4)\to \tc \to 0
$$
and \ref{E_Pi} below which says that $h^0(E_{\Pi}(-e-4))=0$ implies $k\geq 0$. Hence we find
$\beta^{**}: \tc ^{**}\simeq \oc_{\Pi}(k-2)\hookrightarrow\oc_{\Pi}(-2-l)$ with $l\geq 0$ and $k\geq 0$
hence $l=k=0$ and $\beta^{**}$ is an isomorphism.
\end{proof}

By \ref{l=0}, $\Sc$ is the ideal sheaf of a point $q\in \Pi$ and $\beta $ is an isomorphism in $\tilde{\Pi}:= \Pi - (\Phi\cup q)$. 

We can now prove (modulo Lemma \ref{E_Pi} below) the main result of this section.

\begin{proposition}
\label{e=1}
With notations as above we have $e=1$
\end{proposition}

\begin{proof}
Since $\beta $ is an isomorphism on $\tilde{\Pi}$ then 
$\alpha: \nc^*_Z\mid_{\tilde{\Pi}}\to \nc^*_{\Pi}\mid_{\tilde{\Pi}}$ has rank one at any point. This shows that $Z$ is \textit{curvilinear} at any point of $\tilde{\Pi}$. Consider a general $\Pt\subset \Pq$ (not meeting $\Phi\cup q$). Then the curve
$Y=Z\cap \Pt$ is a primitive multiple structure on a line $L$. Furthermore, since $c_1(\Sc)=0$
then the double structure contained in $Y$ is given, via Ferrand's construction, by a surjection  $\nc^*_L\to \oc_L$. Now let $\Gamma$ be a primitive multiple structure, of multiplicity $k$, on a smooth curve $C$. If $I$ is the ideal defining $C$ in $\Gamma$, then the graded ring associated to the $I$-adic filtration of $\oc _{\Gamma}$ is $\bigoplus _{i=0}^k \lc ^{\otimes i}$ ($\lc = I/I^2$) and $\lc ^{\otimes (k-1)} = Hom_{\oc _{\Gamma}}(\oc _C,\oc _{\Gamma})$, so $\lc ^{\otimes (k-1)} \simeq \omega _C\otimes \omega _{\Gamma}^{-1}$. It follows that in our case: 
$\oc _L \simeq \omega _L(e+1)$ and we are done.
\end{proof}

To conclude this section let's prove that $h^0(E_{\Pi}(-e-4))=0$.\\
If $H$ is an hyperplane through $\Pi$, then 
$h^0(E_H(-e-2))\not = 0$ and by \ref{E_H} a section $\sigma _H\in h^0(E_H(-e-2))$ must vanish in codimension two:
$$ 0 \to \oc_H \stackrel{\overline{\sigma}_H}{\to} E_H(-e-2) \to \ic _{X_H}(-e+1) \to 0$$
where  $X_H$ is a curve supported in $\Pi $.

Now we restrict $\overline{\sigma}_H$ to $\Pi$, again this section vanishes on a divisor, after division by the equation of this divisor, we get a section vanishing in codimension two:
$$ 0 \to \oc_{\Pi} \stackrel{\overline{\sigma}_{\Pi}}{\to} E_{\Pi}(-e-2-t) \to \ic _{\Delta}(-e+1-2t) \to 0 \:\: (*)$$
where  $t \geq 1$.

\begin{lemma}
\label{E_Pi}
We have $h^0(E_{\Pi}(-e-4))=0$ hence $t=1$. 
\end{lemma}
\begin{proof} We distinguish two cases: (a) $dim(S \cap \Pi )=0$, (b) $dim(S \cap \Pi )=1$.\\
(a) Since $-e+1-2t <0$, from $(*)$ $h^0(E_{\Pi}(-e-2-t))=1$ and this is the least twist of $E_{\Pi}$ having a section. By restricting the section $s$ to $\Pi$, we get: $h^0(\ic _{S\cap \Pi}(3-t))\neq 0$. If  $t \geq 2$  it follows that $t=2$ and $h^0(\ic _{S \cap \Pi}(1)) \neq 0$, but this is absurd ($S$ has no $d$-secant line).\\
(b) Let $P$ denote the one-dimensional part of $S \cap \Pi$. Restricting $s$ to $\Pi$, after division by the equation of $P$, we get:
$$ 0 \to \oc_{\Pi} \stackrel{\overline{s}_{\Pi}}{\to} E_{\Pi}(-p) \to \ic _{\Gamma}(e+5-2p) \to 0 $$
where $\Gamma$ is zero-dimensional. Twisting by $p-e-2-t$ and since $h^0(\ic _{\Gamma}(3-2p-t))=0$, from $(*)$ we get: $p=e+2+t$. In particular $p \geq e+3$. The next lemma shows that $p=e+3$, hence $t=1$. \end{proof}

\begin{lemma}
\label{twotwo}
Let $S \subset \Pq$ be a smooth surface with $\omega _S \simeq \oc _S(e)$. If $S$ contains a plane curve, $P$, then $deg(P) \leq e+3$. Moreover if $deg(P)=e+3$, then $P^2=0$.
\end{lemma}

\begin{proof} See \cite{F}. For convenience of the reader we give a proof. Let $H$ be an hyperplane containing $P$, we have $C=H \cap S=P \cup Y$ and we may assume that no irreducible component of $Y$ is contained in $P$ (see \cite{F}). We have $\omega _{C|P} = \omega _P(\Psi) = \oc _P(p-3+\Psi)$ ($\Psi =P \cap Y$, $p=d(P)$). Since $\omega _C\simeq \oc _C(e+1)$, $\oc _P(p-e+\Psi ) \simeq \oc _P(e+1)$ and we conclude because $\Psi > 0$.\\
The last statement follows by adjunction. 
\end{proof}

Let's note in passing, for later use, the following consequence:

\begin{corollary}
\label{S.Pi=1}
If $Z$ contains a plane $\Pi$ such that $dim(S \cap \Pi )=1$ and if $e \geq 1$, then $\Pi \subset Sing(\Sigma )$.
\end{corollary}

\begin{proof} If $\Pi$ has equations $x_0=x_1=0$, then $\Sigma$ has an equation of the form $F=x_0f+x_1g=0$ where $f,g$ are cubic forms. The partials $F_i$ vanish on $\Pi$ if $i>1$ and $F_0|\Pi =f|\Pi$, $F_1|\Pi =g|\Pi$. Since $P \subset S \cap Z$, $P \subset Sing(\Sigma )$, so $P \subset (f \cap g)|\Pi$, for degree reasons ($deg(P)=e+3$ by  \ref{E_Pi}) this implies $f|\Pi=g|\Pi=0$ and we are done. 
\end{proof}


\section{Double structures on degree three surfaces.}
\label{dble-str-surfaces}

\begin{lemma}
\label{dble_scroll}
Let $Z \subset \Pq$ be a locally complete intersection double structure on a cubic scroll $T$. If $\omega _Z \simeq \oc _Z(-\alpha)$, then $\alpha =0$.
\end{lemma}

\begin{proof}
The double structure corresponds to a quotient: $N^*_T \to \lc \to 0$. It follows that $\lc ^*$ is a subbundle of $N_T$, hence $c_2(N_T\otimes \lc )=0$. If $\lc = \oc _T(aC_0+bf)$ (notations as in \cite{Ha}), then we get:
$$b= \frac{a^2-4a-9}{2a+3}\:\: (*)$$
Since $Z$ is l.c.i., we have, by local algebraic linkage, $\lc = \ic _{T,Z} = \omega _T \otimes \omega _Z^{-1}$, hence $\lc = \omega _T(\alpha )$. It follows that $a = \alpha -2$ and $b= 2\alpha -3$. Plugging into $(*)$ we get the result.
\end{proof}

Now let $T \subset \Pq$ denote a cone over a twisted cubic. As before, if $Z$ is a l.c.i. double structure on $T$, then $Z$ is defined by the Ferrand construction (observe that $\ic _T^2 \subset \ic _Z$)  and we have: $N^*_T \to \lc \to 0$ where $\lc \simeq \ic _{T,Z}$, moreover, if $Z$ is subcanonical, we have: $\lc \simeq \omega _T(-\alpha )$.
\par
Denote by $p: \tilde{T} \to T$ the desingularization of $T$.

\begin{lemma}
\label{Ncub_cone}
With notations as above:
\begin{enumerate} \item $p^*(N^*_T) \simeq \omega _{\tilde{T}} \oplus \omega _{\tilde{T}}$ (modulo torsion).
\item $p^*(\omega _T) \simeq \omega _{\tilde{T}}$ (modulo torsion).
\end{enumerate}
\end{lemma}

\begin{proof}
 1) This follows from the fact that if $C \subset \Pt$ is a twisted cubic, then $N^*_C \simeq \oc _{\Pu}(-5) \oplus \oc _{\Pu}(-5) \simeq \omega _C(-1) \oplus \omega _C(-1)$
\par \noindent
2) This can be seen arguing as above, here is another proof: the cone $T$ is linked to a plane $P$ by a complete intersection, $U$, of type $(2,2)$ and we have an exact sequence: $0 \to \ic _U \to \ic _P \to \omega _T(1) \to 0$. Now, $\omega _T(1) \simeq \oc _T(-D)$, where $D$ is the curve $P \cap T$ which is clearly the union of two rulings of $T$. Now (with notations as in \cite{Ha}): $p^*(-D-H) = -2C_0 -5f = K_{\tilde{T}}$, which shows that $p^*(\omega _T)/torsion \simeq \omega _{\tilde{T}}$.
\end{proof}

\begin{lemma}
\label{dble_str_cubcone}
If $Z$ is a l.ci., subcanonical double structure on the cubic cone $T$, then $Z$ is a complete intersection of type $(2,3)$ and $\omega _Z \simeq \oc _Z$.
\end{lemma}

\begin{proof}
 The morphism $N^*_T \to \lc \simeq \omega _T(-\alpha ) \to 0$, yields: $p^*(N^*_T) \to p^*(\omega _T(-\alpha ) \to 0$, killing the torsion this gives by Lemma \ref{Ncub_cone}: $0 \to \fc \to \omega _{\tilde{T}} \oplus \omega _{\tilde{T}} \to \omega _{\tilde{T}} (-\alpha ) \to 0$, where $\fc$ is locally free. Looking at determinants, we conclude that $\fc \simeq \omega _{\tilde{T}}(\alpha)$. It follows that $\alpha =0$. Now taking a general hyperplane section, we get a double structure, $X$, on a twisted cubic, $C$, such that: $\ic _{C,X} \simeq \omega _C(-1) \simeq \oc _{\Pu}(-5)$. It is easily seen that such a double structure is a complete intersection (observe that $h^0(\ic _X(2))=1$ and then show there is a cubic , not multiple of the quadric, containing $X$).
\end{proof}


\section{End of the proof.}
\label{epf}

At this point, by  \ref{e=1}, we may assume (if $e\geq 2$ and if $S$ is not a complete intersection) that $Z_{red}$ contains an irreducible, non degenerate, cubic surface, $T$. Observe that, since $deg(Z) \leq 9$ and since every irreducible component of $Z_{red}$ appears with multiplicity in $Z$ (Lemma \ref{gen_Z}), if $T \subset Z_{red}$, we have the following possibilities:
\begin{enumerate}
\item $Z_{red}=T$ and $Z$ is a double or triple structure on $T$.
\item $Z$ is a double structure on $T$ union a double (or triple) structure on a plane $\Pi$.
\end{enumerate}

\begin{lemma}
\label{triple-cubic-smooth}
With notations as above, if $Z_{red}=T$ and $d(Z)=9$, then $T$ is a cubic cone.
\end{lemma}

\begin{proof}
If not $T$ is a cubic scroll. Since $Z$ is a l.c.i. triple structure on a smooth support, it is a primitive structure, hence contains a double structure, $Z_2$, given by $\nc ^*_T \to \lc \to 0$, where $\lc$ is a line bundle on $T$. By local algebraic linkage we have $\lc ^{\otimes 2} \simeq \omega _T(e+2)$. Taking as basis of $Pic(T)$ the ruling $f$ and the hyperplane class $h$, we have $\omega _T \sim f-2h$, hence $\omega _T(e+2) \sim f+eh$ which is not divisible by two.
\end{proof}

\begin{lemma}
\label{notriple}
With notations as above $Z$ cannot be a triple structure on $T$.
\end{lemma}

\begin{proof}
If $Z$ is a triple structure on $T$, then $T$ is a cubic cone (see \ref{triple-cubic-smooth}). By \ref{inSing} $T$ is contained in the singular locus of $\Sigma$. So if $F =\Sigma \cap H$ is a general hyperplane section, $F$ is a quartic surface with a twisted cubic $\Gamma = T\cap H$ in its singular locus. Such a quartic surface is ruled. Let $v$ denote the vertex of $T$, we claim that $\Sigma$ is a cone of vertex $v$, base $F$. Assume the claim for a while.\\
If $\Pi$ is a general plane lying on $\Sigma$, then $E_{\Pi}(-e-3)$ has a section dividing the restrictions to $\Pi$ of $s$ and $\sigma$ (the sections giving $S$ and $Z$). This follows because the sections are proportional on $\Sigma$ and $\Pi$ intersects $T$ in two lines (not necessarily distinct). This means that $\Pi \cap S$ contains a plane curve, $P$, of degree $e+3$. If we let $\Pi$ vary we get a family of rationally equivalent plane curves inside $S$, with vanishing self intersection by lemma \ref{twotwo}. This implies that $P$ doesn't pass through $v$. Now both $S$ and $Z$ meet every plane of the ruling of $\Sigma$ along a curve (for dimension reason), these two plane curves intersect outside of $v$, moving the plane we get that $dim(S \cap Z) > 0$. We conclude with Lemma \ref{gen_Z}.\\
{\it Proof of the claim:} The quartic $F$ is ruled by bisecants or tangents to $\Gamma$. In the first case, the plane spanned by $v$ and a bisecant meet $\Sigma$ along two double lines through $v$ (recall that $T \subset Sing(\Sigma )$) and the bisecant; so this plane is contained in $\Sigma$. In the second case let $\Sigma '$ denote the developable of $T$, then $\Sigma$ and $\Sigma '$ have the same intersection with $H$, since $H$ is general, we get $\Sigma = \Sigma '$.
\end{proof}

\begin{proof} (of Theorem \ref{thm_subcan_quartic})\\
As observed in \ref{gen_setting} we may assume $e \geq 2$, so to prove the theorem it is enough to show that case 1. and case 2. above are impossible.\\
{\em Case (1):} By Lemma \ref{notriple} we may assume that $Z$ is a double structure on $T$. If $T$ is a cone, we conclude with Lemma \ref{dble_str_cubcone}; if $T$ is a scroll, with Lemma \ref{dble_scroll}.\\
{\em Case (2):} By Lemma \ref{inSing}, $T \subset Sing(\Sigma )$, so $\Sigma$ is ruled in planes. Since $dim(T \cap \Pi)=1$ (otherwise $Z$ is not locally Cohen-Macaulay, \cite{Ha3}), $\Pi$ is a plane of the ruling. This implies $dim(S \cap \Pi)=1$. By Corollary \ref{S.Pi=1}, $\Pi \subset Sing(\Sigma )$. So $Z_{red} \subset Sing(\Sigma )$, but this is impossible because $deg(Z_{red})=4$. 
\end{proof}

\bigskip



\begin{thebibliography}{subsurfPq}
\bibitem{BC} Ballico, E.-Chiantini, L.: {\it On smooth subcanonical varieties of codimension two in $\Pn$, $n \geq 4$}, Annali Mat. Pura Appl., 135, 99-118 (1983) 
\bibitem{EF} Ellia, Ph.-Franco, D.: {\it On codimension two subvarieties of $\Pq$, $\Pcq$}, J. Algebraic Geometry, 11, 513-533 (2002)
\bibitem{E} Ellia, Ph.: {\it A survey on the classification of codimension two subvarieties of $\Pn$}, Le Matematiche, vol. LV, 285-303 (2000)
\bibitem{EP} Ellingsrud, G.-Peskine, Ch.: {\it Sur les surfaces lisses de $\Pq$}, Invent. Math., 95, 1-11 (1989)
\bibitem{F} Folegatti, C.: {\it On a special class of smooth codimension two subvarieties in $\Pn$, $n \geq 5$}, preprint (2003)
\bibitem{Ha} Hartshorne, R.: {\it Algebraic geometry}, G.T.M., 52, Springer (1977)
\bibitem{Ha3} Hartshorne, R.: {\it Complete intersections and connectedness}, Amer. J. of Math., 84, 497-508 (1962) 
\bibitem{K} Koelblen, L.: {\it Surfaces de $\Pq$ trac\'ees sur une hypersurface cubique} J. reine u. angew. Math., 433, 113-141 (1992)

\end{thebibliography}
\end{document}